\documentclass[a4paper,12pt]{article}

\usepackage{amsmath}
\usepackage{amsthm}
\usepackage{amssymb}  
\usepackage{latexsym} 
\usepackage[all]{xy}
\usepackage{comment}
\usepackage{rotating}
\usepackage{slashbox}
\usepackage{multirow}
\usepackage{rotating}

\newcommand{\pp}{{\mathfrak {p}}}

\newcommand{\CC}{{\mathcal C}}

\newcommand{\mm}{\mathfrak m}

\newcommand{\Char}{\operatorname{char}}

\newcommand{\Spec}{\operatorname{Spec}}

\newcommand{\mult}{\operatorname{mult}}

\newcommand{\GL}{\operatorname{GL}}

\newcommand{\Kbar}{\overline{K}}
\newcommand{\kbar}{\overline{k}}

\newcommand{\isom}{ \cong }

\newcommand{\OK}{{\mathcal{O}_K}}

\newcommand{\Gal}{\operatorname{Gal}}

\newcommand{\Frac}{\operatorname{Frac}}

\newcommand{\PP}{{\mathbb P}}
\newcommand{\Q}{{\mathbb Q}}

\newcommand{\Hom}{\operatorname{Hom}}

\newcommand{\Z}{{\mathbb Z}}

\newfont{\wncyr}{wncyr10 at 12pt}
\newfont{\wncyrten}{wncyr10 at 10pt}

\newenvironment{Proof}{\par\noindent{\sc Proof:}}%
                      {\hspace*{\fill}\nobreak$\Box$\par\medskip}
                       {\hspace*{\fill}\nobreak$\Box$\par\medskip}
\newenvironment{myitemize}
{\begin{itemize}
\setlength{\itemsep}{1pt}
\setlength{\parskip}{0pt}
\setlength{\parsep}{0pt}}
{\end{itemize}}

\newtheorem{Proposition}{Proposition}[section]
\newtheorem{Theorem}[Proposition]{Theorem}
\newtheorem{Lemma}[Proposition]{Lemma}
\newtheorem{Corollary}[Proposition]{Corollary}

\theoremstyle{definition}

\newcounter{nootje}
\begin{document}

\title{On Quadratic Twists of Hyperelliptic Curves}
\author{Mohammad Sadek}
\date{}
\maketitle
\begin{abstract}{\footnotesize
Let $C$ be a hyperelliptic curve of good reduction defined over a discrete valuation field $K$ with algebraically closed residue field $k$. Assume moreover that $\Char k\ne2$. Given $d\in K^*\setminus K^{*2}$, we introduce an explicit description of the minimal regular model of the quadratic twist of $C$ by $d.$ As an application, we show that if $C/\Q$ is a nonsingular hyperelliptic curve given by $y^2=f(x)$ with $f$ an irreducible polynomial, there exists a positive density family of prime quadratic twists of $C$ which are not everywhere locally soluble.}
\end{abstract}

\section{Introduction}
\label{sec:intro}

Let $C$ be a nonsingular hyperelliptic curve defined over $\Q$ with an affine model given by the equation $y^2=f(x)\in\Z[x]$, where $\deg(f)\ge3$. The genus of $C$ will be called $g$. If $d>1$ is a square free integer, then we write $C_d$ for the quadratic twist of $C$ by $d$. In particular, $C_d$ is defined by $dy^2=f(x)$.

We try to find an explicit description of the minimal regular model of $C_d/\Q_p$ in terms of the minimal regular model of $C/\Q_p$ itself, when $C$ has good reduction.

Let $\Delta$ denote the minimal discriminant of $C/\Q_p$, see \S 2 of \cite{LiuModeles}. For every prime $p\nmid \Delta$, $C/\Q_p$ has good reduction. Hence the minimal regular model of $C/\Q_p$ is smooth over $\Z_p$. The minimal regular model of $C_d/\Q_p,\;p\mid d,$ is obtained as the minimal desingularisation of a quotient of a smooth scheme by a twisted action of some finite group, see \S \ref{sec:minimal models}. In fact, we prefer to handle the problem over the maximal unramified extension of $\Q_p$ to avoid any complications which might appear because the residue field $\mathbb{F}_p$ is not algebraically closed.

As we investigate minimal regular models of quadratic twists of a hyperelliptic curve $C$ over $\Q_p$, we cannot see how to do the curves $C_d$, where one of the prime divisors $p$ of $d$ is a bad prime of $C$ and $f$ has no simple root when reduced modulo $p$. The difficulty lies in the wide range of possibilities of the structure of the minimal regular model of $C/\Q_p$ when $C$ has bad reduction. Furthermore, we are not aware of any reference which discusses the desingularization of quotient singularities of models of algebraic curves when these singularities are not ordinary double points.

 Now assume $C$ is a hyperelliptic curve defined over $\Q$ given by the equation $y^2=f(x)$, where $f(x)$ is an irreducible polynomial. Using the description of the minimal regular model of a quadratic twist of $C/\Q_p$, we show that there is an infinite number of quadratic twists of $C$ with no $\Q_p$-rational points. In particular, for a nonsingular hyperelliptic curve $C/\Q$, there exists an infinite number of quadratic twists $C_d$ of $C$ such that $C_d(\Q_p)=\emptyset$ for some prime $p$, and hence $C_d(\Q)=\emptyset$.

\section{Hyperelliptic Curves}
\label{sec:hyper}
The material in this section can be found in \S7.4.3 of \cite{Liubook}.

We assume $K$ is a field with $\Char K\ne2$ and algebraic closure $\Kbar$. Two hyperelliptic equations with coefficients in $K$
\[y^2=f(x)\textrm{ and }z^2=f'(u)\]
represent isomorphic curves if and only if
\[x=\frac{au+b}{cu+d},\;y=\frac{ez}{(cu+d)^{g+1}}\textrm{ where }\left(\begin{array}{cc}a & b \\c & d \\\end{array}\right)\in \GL_2(K),\;e\in K^*.\]
We associate a discriminant $\Delta$ to a hyperelliptic equation $y^2=f(x)$ as in \S 2 of \cite{LiuModeles}. This equation defines a smooth curve if and only if $\Delta\ne0$.

 By a {\em hyperelliptic curve} $C$ over $K$ we mean a smooth curve of genus $g\ge2$ endowed with a morphism $C\to\PP_K^1$ of degree $2$. There exists a hyperelliptic equation $y^2=f(x)\in K[x]$ describing $C$ with $\deg f\in\{2g+1,2g+2\}$. The fact that $C$ is smooth implies $f(x)$ is separable over $\overline{K}.$ The equation $y^2=f(x)$ has one singularity at infinity. If $\deg f=2g+1$, the singularity at infinity corresponds to one point $\infty$ on the hyperelliptic equation. If $\deg f=2g+2$, the singularity corresponds to two points $\infty^+$ and $\infty^-$ on the hyperelliptic equation, and these can be distinguished by the value of the rational function $y/x^{g+1}.$ If $\deg f=2g+2$, then there exists a hyperelliptic equation $z^2=f'(u)$ describing $C$ over $\overline{K}$ with $\deg f'=2g+1$, and it describes $C/K$ if and only if $f(x)$ has a zero in $K$.

Let $C$ be a smooth hyperelliptic curve of genus $g\ge2$ defined over $K$ by the equation $y^2=f(x)$. We will denote the hyperelliptic involution on $C$ by $\iota:C\to C$. Let $K' = K(\sqrt{d})$ be a separable quadratic extension of $K$ where $d\in K\setminus K^{* 2}$. By a quadratic twist $C_d$ of $C$ we mean the hyperelliptic curve obtained from the curve $C/K'$ by twisting the curve $C/K$
by the cohomology class corresponding to $K'$ in $H^1(K,\langle\iota\rangle).$

This means that if $\sigma$ generates $\Gal(K'/K)$, then the twisted action of $\Gal(K'/K)$ on $C(K')$ is given by $Q\to \iota(\sigma Q)$. To produce an explicit equation describing $C_d$, we consider the quadratic character $\chi:\Gal(K'/K)\to\{\pm1\}$ associated with $K'/K$, i.e., $\chi(\sigma)=\sqrt{d}^{\sigma}/\sqrt{d}$. Then we define a cocycle in
\[H^1(K,\langle\iota\rangle)\isom H^1(K,\Z/2\Z)=\Hom(\Gal(\overline{K}/K),\Z/2\Z)=\Hom(\Gal(\overline{K}/K),\pm1)\]
by \[\xi:\Gal(\overline{K}/K)\to \{\pm1\};\;\xi_{\tau}=[\chi(\tau)].\]
Now, $\overline{K}(C)=\overline{K}(x,y)$ and $\overline{K}(C_d)=\overline{K}(x,y)_{\xi}$ (the twist of the function field by the cocycle $\xi$). Since $\iota(x,y)=(x,-y)$, the action of $\sigma$ on $\overline{K}(x,y)_{\xi}$ is described by
\[\sqrt{d}^{\sigma}=\chi(\sigma)\sqrt{d},\;x^{\sigma}=x,\;y^{\sigma}=\chi(\sigma)y.\]
Thus the functions which are fixed by $\Gal(K'/K)$ in $\overline{K}(x,y)_{\xi}$ are $x'=x,\; y'=y/\sqrt{d}$, hence they are in $K(C_d)$. They satisfy the equation
\[dy'^2=f(x').\]
The curves $C$ and $C_d$ are isomorphic over $K'$ via $(x',y')\mapsto(x',y'\sqrt{d})$.

\section{Minimal regular models of hyperelliptic curves}
\label{sec:minimal models}

We assume that $K$ is a complete discrete valuation field with ring of integers $\OK$, valuation $\nu$, uniformiser $t$ and  residue field $k$ with $\Char k\ne 2$. Set $S=\Spec \OK$.

Let $C$ be a hyperelliptic curve defined over $K$. In \cite{LiuModeles}, Liu associates to $C$ a projective model $W$, a {\em Weierstrass model} of $C$, defined over $\OK$, arising from a hyperelliptic equation of $C$ with integral coefficients. The discriminant $\Delta_{W}$ of $W$ is the discriminant of this hyperelliptic equation. The model $W$ is said to be {\em minimal} if $\Delta_{W}$ is minimal, i.e., $\nu(\Delta_{W})$ is the least possible valuation among the valuations of the discriminants of the hyperelliptic equations related to our equation via the transformations of the form given in \S \ref{sec:hyper}. If $C(K)\ne\emptyset$, then $W$ being minimal implies that the minimal regular model of $C$ is the minimal desingularisation of $W$, see Corollaire 5 of \cite{LiuModeles}.

 Let $W$ be a minimal Weierstrass model of $C$. The curve $C$ has good reduction if $\nu(\Delta_{W})=0.$ In fact, the latter statement is equivalent to saying that the minimal regular model $\CC$ of $C$ over $S$ is smooth, see \S 3 of \cite{LiuModeles}. Moreover, $\CC$ is the unique smooth model of $C$ over $S$, (\cite{Liubook}, Proposition 10.1.21 (b)).

\begin{Lemma}
\label{lem2}
 Assume that $C$ has good reduction over $k$. Let $K'/K$ be a finite extension with residue extension $k'/k$. Then $C\times_KK'$ has good reduction over $k'$.
\end{Lemma}
\begin{Proof}
Let $\CC$ be a minimal Weierstrass model of $C$. Let $\nu'$ be the valuation corresponding to $K'$. Since $\nu'=e_{K'/K}\times\nu$, where $e_{K'/K}$ is the ramification index of $K'/K$, one has $\nu'(\Delta_{\CC})=e_{K'/K}\times\nu(\Delta_{\CC})=0$.
\end{Proof}

Recall the following results, which allow us to determine whether a curve
defined over a complete discrete valuation field $K$ has a $K$-rational point.

\begin{Lemma}[\cite{Liubook}, Corollary 9.1.32]
\label{lem1}
Let $C$ be an algebraic curve of genus $g\ge1$ defined over $K$. Let $\CC\to S$ be the minimal regular model of $C$. Assume that $C(K)\ne\emptyset$. Then a point $P\in C(K)$ is reduced to a point $\widetilde{P}\in \CC_k(k)$, and $\CC_k$ is smooth at $\widetilde{P}$. In particular, $\widetilde{P}$ belongs to a single irreducible component of multiplicity $1$ in $\CC_k$.
\end{Lemma}
\begin{Lemma}
\label{lem3}
Let $C$ be a smooth hyperelliptic curve over $K$. Assume that $C$ has good reduction over $k$. Let $K'/K$ be a quadratic extension with residue field $k'/k$. Let $\CC'$ be the minimal regular model of $C\times K'$. Then $\CC'$ is smooth over $\OK$. Moreover, $\CC'_{k'}$ consists of one irreducible component of multiplicity $1$.
\end{Lemma}
\begin{Proof}
Since $C$ has good reduction over $k$, then it extends to a smooth relative curve $\CC/_{\mathcal{O}_K}$. This relative curve is the minimal regular model of $C$ over $\mathcal{O}_{K}$. Hence it consists of one irreducible component $\Gamma$. If $C(K)\ne\emptyset$, then $\mult(\Gamma)=1$, see Lemma \ref{lem1}, otherwise $\mult(\Gamma)=2$. The reason for the latter statement is as follows: Let $D$ be the image of $\Gamma$ in $\PP_S^1$ under the morphism $g:\CC\to\PP^1_S$. We will denote the generic points of $\Gamma$ and $D$ by $\xi$ and $\xi'$ respectively. The morphism $g$ restricts to $\mathcal{O}_{\Gamma,\xi}\to\mathcal{O}_{D,\xi'}$. The valuations $\nu_{\Gamma}$ and $\nu_{D}$ are the corresponding normalised valuations to $\Gamma$ and $D$ respectively. Remember that $\nu_D(t)=1$ because $D$ is reduced. One has $[\Frac(\mathcal{O}_{\Gamma,\xi}):\Frac(\mathcal{O}_{D,\xi'})]=2$. Since $\Gamma$ is not reduced, we deduce that $t$ ramifies in $\Frac(\mathcal{O}_{\Gamma,\xi})$. Thus $\nu_{\Gamma}(t)=2$.

Lemma \ref{lem2} implies that the minimal regular model of $C\times K'$ is smooth. Again it consists of one irreducible component $\Gamma'$ of multiplicity $1$. This is clear if $\mult(\Gamma)=1$. If $\mult(\Gamma)=2,$ then $\mult(\Gamma')=\mult(\Gamma)/[K':K]$, see for example \S 2.4 of \cite{LorenziniModels}. It is true that the mentioned reference gives result when $\kbar=k$, but if we take a base change over the maximal unramified extension of $K$, then the multiplicity of components will not change.
\end{Proof}

We have to mention that the statements of Lemmas \ref{lem2} and \ref{lem3} are true for any base extension. In other words, smoothness is preserved by arbitrary base change. We wrote down the proofs when the base change is quadratic for the convenience of the reader.

In what follows we assume $k$ is algebraically closed. Hence $K'=K(\sqrt{t})$ is the unique quadratic extension of $K$, and it is totally and tamely ramified. Furthermore, $K'/K$ is Galois. Let $G:=\Gal(K'/K)=\langle\sigma\rangle.$

 Again $C/K$ is a hyperelliptic curve. We assume that $C/K$ has good reduction. We are concerned with the twisted action of $\sigma$ on $C(K')$ given by $Q\to \iota(\sigma Q)$, where $\iota:C\to C$ is the hyperelliptic involution on $C$ and $\sigma(Q)$ is the usual Galois action of $\sigma$ on $C(K')$. Now the automorphism $\sigma:C\to C$ extends to the minimal regular models $\CC$ and $\CC'$ of $C$ and $C\times K'$ respectively. We will denote the extended automorphism by $\sigma$ again.

Let $C_t$ be the quadratic twist of $C$ by $t$. In what follows we obtain the minimal regular model of $C_t/K$ as the minimal  desingularisation of the quotient scheme $\CC'/\langle\sigma\rangle$ by the twisted action of $\sigma$. In other words, we construct the minimal regular model of $C_t/K$ from the minimal regular model of $C/K$.

The first step is to find the fixed points of the twisted action of $\sigma$ on $\CC'$. This is because the singular points of $\CC'/\langle\sigma\rangle$ lie among the images of the points of $\CC'$ fixed by $\sigma$. Let $y^2=f(x)$ be a minimal hyperelliptic equation describing $C$. Then the
assumption that $C$ has good reduction implies that $f(x)$ has no repeated roots over $k$. If $P\in C\times K'$, we denote its reduction by $\widetilde{P}.$


\begin{Proposition}
\label{prop1}
 Let $P\in C\times K'$. The following are equivalent:
 \begin{myitemize}
  \item[(i)] $P$ is fixed under the twisted action of $\sigma$.
  \item[(ii)] $\widetilde{P}=(\tilde{x},\tilde{y})$ is fixed under the twisted action of $\sigma$.
  \item[(iii)] $f(\tilde{x})=0$ over $k$.
  \end{myitemize}
  \end{Proposition}
\begin{Proof}
The twisted action of $\sigma$ on $P$ is given by $\iota(\sigma P)$. Let $P=(x_0+x_1\sqrt{t},y_0+y_1\sqrt{t})$, where $x_0,x_1,y_0,y_1\in K$ is fixed if and only if $(x_0-x_1\sqrt{t},-y_0+y_1\sqrt{t})=(x_0+x_1\sqrt{t},y_0+y_1\sqrt{t}).$ Hence $x_1=y_0=0$. Whence $P$ is fixed if and only if $P=(x_0,y_1\sqrt{t})$ where $x_0,y_1\in K$. The latter is equivalent to $\widetilde{P}=(\tilde{x}_0,0),$ or equivalently $f(\tilde{x}_0)=0$ in $k$. Now if $(x,y)\in \CC_k$, then it is $\sigma$-fixed if and only if $(x,y)=(x,-y)$, i.e., $2y=0$ but $2\in k^*$, hence $y=0.$ So $(i)\Leftrightarrow(ii)$ holds.
\end{Proof}

Recall that $C$ and $C_t$ are isomorphic over $K'$. Therefore both $C\times K'$ and $C_t\times K'$ have the same minimal regular model $\CC'\to\Spec\mathcal{O}_{K'}$. The model $\CC'$ is smooth.
Since $\CC'$ is projective, the $S$-quotient scheme $Z:=\CC'/ \langle\sigma\rangle$ is constructed by glueing together the rings of invariants of $\langle\sigma\rangle$-invariant affine open sets of $\CC'$. Moreover, $Z$ is a normal scheme and hence its singular points are closed points of the special fiber.

\begin{Proposition}
\label{prop2}
Let $\sigma_k:{\CC'}_k\to{\CC'}_k$ and $\sigma_k^{red}:{\CC'}_k^{red}\to{\CC'}_k^{red}$ be the natural morphisms induced by $\sigma$. Then the natural map
\[{\CC'}_k^{red}/\langle\sigma_k^{red}\rangle\to Z_k^{red}:=(\CC'/\langle\sigma\rangle)_k^{red}\]
is an isomorphism over $k$.
\end{Proposition}
\begin{Proof}
See (\cite{Lorenzinidual}, Facts II, III).
\end{Proof}

In fact if $\alpha:\CC'\to Z$ is the quotient map, then $\alpha$ induces a natural map $\CC'_k\to Z_k^{red}$ which factors as follows:
\[\CC'_k\to \CC'_k/\langle\sigma_k\rangle\to Z_k^{red}\]
where the second map is the normalisation map of $Z_k^{red}$, see (\cite{Lorenziniwildquotient}, p.21).

\begin{Proposition}
\label{prop3}
The generic fiber of $Z:=\CC'/\langle\sigma\rangle$ is isomorphic to $C_t/K.$
\end{Proposition}
\begin{Proof}
The generic fiber of $Z$ is given by
\[Z_{K}=Z\times_{\mathcal{O}_K}K=(\CC'/\langle\sigma\rangle)\times_{\mathcal{O}_K}K=(\CC'\times_{\mathcal{O}_{K}}K)/\langle\sigma\rangle.\]
But one has
\[\CC'\times_{\mathcal{O}_K}K=\CC'\times_{\mathcal{O}_{K'}}\mathcal{O}_{K'}\times_{\mathcal{O}_{K}}K=
\CC'\times_{\mathcal{O}_{K'}}K'=C\times K'.\]
Since we consider the twisted action of $\sigma$, we have $Z_{K}=C_t.$
\end{Proof}

Now we aim to prove that the minimal desingularisation $\widetilde{Z}$ of the $S$-quotient scheme $Z=\CC'/\langle\sigma\rangle$ is the minimal regular model of $C_t/K$. Moreover, we will show that this model consists of one irreducible component of multiplicity $2$, the image of the special fiber $\CC'_k$ of $\CC'$ in $Z$, and a finite number of multiplicity-$1$ irreducible components each of which corresponds to a singular point of $Z$. In particular, if $\rho:\CC'\to Z$ is the quotient map, we will prove that since $\CC'$ is smooth then the only singular points of $Z$ are the images of the $\sigma$-fixed points of $\CC'$ under $\rho$.

Consider the morphism
\[\CC'_{k}\xrightarrow{\beta} \CC'_k/\langle\sigma_{k}\rangle.\]
Let $y^2=f(x)$ be a minimal hyperelliptic equation defining $C$. Since $C$ is smooth, the minimal regular models $\CC$ and $\CC'$ are smooth, and they are still defined by this affine equation.  The ramification points of $\beta$ are the zeros of $f(x)$ over $k$ plus the point at infinity when $\deg f=2g+1$. Since $f(x)$ splits completely into linear factors over $k$, because $k=\overline{k}$, the number of the ramification points of $\beta$ is $2g+2$.

\begin{Theorem}
\label{thm1}
Let $C/K$ be a hyperelliptic curve of genus $g$. Assume that $C$ has good reduction over $k$. Let $\CC'$ be the minimal regular model of $C\times K'$. Again $\Gal(K'/K)=\langle\sigma\rangle$. Let $Z$ denote the quotient of $\CC'$ by the twisted action of $\sigma$ by the hyperelliptic involution. Then $Z$ is singular exactly at the images $x_1,\ldots,x_m$, $m=2g+2$, of the ramification points of the morphism $\beta: \CC'_k\to \CC'_k/\langle\sigma_k\rangle$ in $Z$.

Let $\widetilde{Z}\to Z$ be the minimal desingularisation of $Z$. Then $\widetilde{Z}$ is the minimal regular model of the quadratic twist $C_t$ of $C$ over $K'$. Moreover, $\widetilde{Z}_k$ consists of an irreducible component $\Theta$ of multiplicity $2$ and components $\Gamma_1,\ldots,\Gamma_m$ of multiplicity $1$, each of which corresponds to blowing-up one of the $x_i$'s, moreover $\Theta$ and $\Gamma_i$'s are of genus zero, see the following figure.
\end{Theorem}
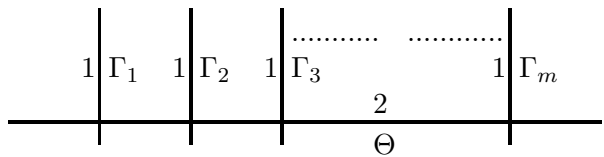
\setlength{\unitlength}{0.6cm}\linethickness{.3mm}
\small
\begin{picture}(18,5)
\put(5,2){\line(1,0){6}}
\put(7,1.5){\line(0,1){3}}\put(6.6,3){$1$}\put(7.2,3){$\Gamma_1$}
\put(9,1.5){\line(0,1){3}}\put(8.6,3){$1$}\put(9.2,3){$\Gamma_2$}
\put(11,1.5){\line(0,1){3}}\put(10.6,3){$1$}\put(11.2,3){$\Gamma_3$}
\put(11.2,3.8){$...........{\hskip10pt}............$}
\put(16,1.5){\line(0,1){3}}\put(15.6,3){$1$}\put(16.2,3){$\Gamma_m$}
\put(11,2){\line(1,0){7}}\put(13,1.3){$\Theta$}\put(13,2.2){$2$}
\end{picture}\\
\begin{Proof}
We have seen that the generic fiber $Z_K$ is the curve $C_t/K$, Proposition \ref{prop3}. Furthermore $Z$ is a normal scheme. As we know that $\CC'_k$ consists of one irreducible component of multiplicity 1, Lemma \ref{lem3}, the image of this irreducible component in $Z_k$ has multiplicity $[K':K]=2$, see Fact IV of \cite{Lorenzinidual}.

To obtain the minimal desingularisation of $Z$ we only need to blow-up the singularities of $Z.$ We will show first that $x_1,\ldots,x_m$ are exactly the singular points of $Z$. If the points $x_1,\ldots,x_m$ were regular, then the morphism $\alpha:\CC'\to Z$ would be flat above $x_1,\ldots,x_m$. Hence by Zariski's Purity Theorem, the branch locus $\alpha(\CC'\setminus U)$, where $U$ is the largest open subscheme of $\CC'$ such that $\alpha|_U:U\to Z$ is \'{e}tale, is of codimension 1, see (\cite{Liubook}, Exercise 8.2.15), a contradiction.

Now we blow up $Z$ at each $x_i,\;i=1,\ldots,m$ to construct the minimal desingularisation $\widetilde{Z}\to Z$. Since $x_i$ is a singular point of $Z$, it is known that if $v_1,v_2$ are local parameters of the local ring $\mathcal{O}_{Z,x_i}$, then the twisted action of $\sigma_k$ on $v_1,v_2$ is described as follows:
\[\sigma_k(v_1)=-v_1,\;\sigma_k(v_2)=-v_2,\]
see (\cite{HuiXue}, Lemma 2.2) or (\cite{Parsin}, Lemma 2). Let $T=\Spec\widehat{\mathcal{O}}_{Z,x_i}$ and $T^G=\Spec\widehat{\mathcal{O}}_{Z,x_i}^G.$ Because $x_i$ is non-regular then $T^G$ is a non-regular local scheme. The monomials $s_1=v_1^2,s_2=v_1v_2,s_3=v_2^2$ are invariant under the action of $\sigma_k$ and any $\sigma_k$-invariant polynomial in $v_1,v_2$ is a polynomial in these. Therefore desingulraising $T^G$ is equivalent to desingularising $\Spec k[s_1,s_2,s_3]/(s_2^2-s_1s_3).$ It is known that we need only one blow-up to desingularise $k[s_1,s_2,s_3]/(s_2^2-s_1s_3)$ at its only singular point corresponding to the maximal ideal $\mm=(s_1,s_2,s_3)$. In fact, the exceptional curve of the blowing-up is one irreducible component $\Gamma_i$ isomorphic to $\PP_k^1$, see for example (\cite{Liubook}, Example 8.1.5). Since $\sigma_k$ acts trivially on $x_i$, the multiplicity of the irreducible component $\Gamma_i$ is $1$, see (\cite{Halletame}, Proposition 7.3 (ii)).

Each $\Gamma_i$ has self-intersection $-2$, whereas the multiplicity-2 component has self-intersection $-g-1\le -3.$ Therefore according to Castelnuovo's criterion $\widetilde{Z}$ contains no exceptional divisors, and hence $\widetilde{Z}$ is the minimal regular model of $C_t/K$.

We conclude by computing the genus of the irreducible components of $\widetilde{Z}_k$. We see that the genus of each of the components $\Gamma_i,\;1\le i\le 2g+2,$ is $0$, because $\Gamma_i\isom\PP_k^1$. We can apply Hurwitz's formula, see (\cite{Liubook}, Remark 10.4.8), in order to find the genus of the multiplicity-$2$ irreducible component $\Theta.$ In fact, if $p_a$ denotes the genus, then we have
\begin{eqnarray}
2p_a(\CC'_k)-2&=&2(2p_a(\Theta)-2)+\sum_{1\le i\le m}1\nonumber\\
&=&4p_a(\Theta)+2g-2\nonumber.
\end{eqnarray}
Therefore we deduce $p_a(\Theta)=0$ and $\Theta\isom\PP_k^1.$
\end{Proof}

In Theorem \ref{thm1} we proved that the quadratic twist $C_t/K$ of a hyperelliptic curve $C/K$ of genus $g\ge2$ and with good reduction has a minimal regular model consisting of a component $\Theta\isom\PP^1_k$ of multiplicity $2$, and $2g+2$ irreducible components $\Gamma_1,\ldots,\Gamma_{2g+2}$ each is isomorphic to $\PP_k^1$ and of multiplicity $1$. The intersection numbers are given by $\Theta.\Gamma_i=1$ for every $i$ and $\Gamma_i.\Gamma_j=0$ for every $i\ne j$. In case $g=1$, this is reduction type ${\rm I_0^*}$, see (\cite{Sil2}, Chapter IV, \S 9). While for $g=2$, this is reduction type $[{\rm I}_{0-0-0}^*]$ given for example in p.155 of \cite{Ueno}.

\section{Quadratic twists which are not everywhere locally soluble}
\label{subsec1}

We start this subsection with the following result on irreducible polynomials over $\Q$.

\begin{Lemma}
\label{lem4}
Let $f(x)\in\Z[x]$ be an irreducible polynomial over $\Q$. Then there exists an infinite set of primes $S_f$ in $\Q$ with positive density such that $f(x)$ has no linear factors over $\mathbb{F}_p$ for every $p\in S_f$.
\end{Lemma}
\begin{Proof}
There exists an infinite set of primes $S_f$ such that $f(x)$ has no linear factors over $\mathbb{F}_p$ for every $p\in S_f$, see (\cite{MilneAlgebraicNumber}, Remark 8.40(d)).
 Chebotarev Density Theorem implies that the density $\delta(S_f)$ of $S_f$ exists and satisfies $\delta(S_f)>0$, see for example Exercise 11.3.7 of \cite{Murty}.
\end{Proof}

 Let $C$ be a hyperelliptic curve over $\Q$ of genus $g\ge1$. Then $C$ has good reduction over all but finitely many finite places of $\Q$, see (\cite{Liubook}, Proposition 10.1.21 (a)). The finite set of primes in $\Q$ of bad reduction of $C$ will be denoted by $S_{\Delta}$.

Let $C_d$ be a quadratic twist of $C$ with $d>1$ a square free integer. Since $C$ and $C_d$ are isomorphic over $K=\Q(\sqrt{d})$, it follows that they have the same minimal regular model over the ring of integers of any completion of $K$ at one of its finite places. If $p\not\in S_{\Delta}$, then
this implies that for every prime $\pp\in K$ lying above $p,$ both $C\times K_{\pp}$ and $C_d\times K_{\pp}$ have good reduction, see Lemma \ref{lem2}, where $K_{\pp}$ denotes the completion of $K$ at $\pp$. In particular, the minimal regular model of $C_d\times K_{\pp}$ is smooth.

We will denote the maximal unramified extension of the $p$-adic field $\Q_p$ by $\Q_p^{un}$. The residue field of $\Q_p^{un}$ is $\overline{\mathbb{F}}_p$. We write $\Z_p^{un}$ for the ring of integers of $\Q_p^{un}$. The following fact follows directly from Theorem \ref{thm1}.

\begin{Corollary}
\label{cor3}
Let $C:y^2=f(x)\in\Z[x]$ be a hyperelliptic curve of genus $g\ge1$ over $\Q$. Assume moreover that $f(x)$ is of even degree and irreducible over $\Q$. Let $C_d$ be a quadratic twist of $C$ with $d>1$ a square free integer. If $d$ has an odd prime factor $p\in S_f\setminus S_{\Delta}$, then $C_d(\Q^{un}_p)=\emptyset$. Hence $C_d(\Q_p)=C_d(\Q)=\emptyset.$

In particular, there is a positive density family of prime quadratic twists of $C$ which are not everywhere locally soluble, and hence have no $\Q$-rational points.
\end{Corollary}
\begin{Proof}
Note that $K:=\Q_p^{un}(\sqrt{d})=\Q_p^{un}(\sqrt{p})$ because the residue field $\overline{\mathbb{F}}_p$ is algebraically closed. Since $p\not\in S_{\Delta}$, one has $C\times K$ is smooth and $f(x)$ factors completely into linear factors over $\overline{\mathbb{F}}_{p}$.


Let $\CC\to\Spec\OK$ be the minimal regular model of $C\times K$ and
$\CC_d\to\Spec\Z^{un}_p$ the minimal regular model of $C_d/\Q^{un}_p$. According to Theorem \ref{thm1}, the special fiber $(\CC_d)_{\overline{\mathbb{F}}_{p}}$ consists of an irreducible component of multiplicity $2$ and multiplicity-$1$ irreducible components $\Gamma_1,\ldots\Gamma_m$, $m=2g+2$, corresponding to blowing-up the singular points of $\CC/\langle\sigma\rangle$, where $\Gal(K/\Q_p^{un})=\langle\sigma\rangle$ and the action of $\sigma$ on $\CC$ is the twisted action by the hyperelliptic involution introduced in \S \ref{sec:minimal models}. These singular points correspond to the simple roots of $f(x)$ over $\overline{\mathbb{F}}_{p}$ (plus the point at infinity if $\deg f=2g+1$). But since $f(x)$ has no simple root over $\mathbb{F}_{p}$ as $p\in S_f$, and $\deg f=2g+2$, it follows that each $\Gamma_i$ is defined over $\overline{\mathbb{F}}_p$ and none of the $\Gamma_i$'s is defined over $\mathbb{F}_p$. Indeed, $\CC_d\to\Spec\Z_p$ consists only of the multiplicity-2 component. According to Lemma \ref{lem1}, $C_d(\Q_p)=\emptyset.$
\end{Proof}

\bibliographystyle{plain}
\footnotesize
\bibliography{thesisreferences}
\end{document}